\documentclass[pt12egno]{article}
\usepackage{amsmath, amssymb}
\usepackage{latexsym}

\usepackage{a4,epsfig}
\usepackage{amsfonts}
\usepackage{graphicx}
\usepackage{amscd}
\usepackage{amsbsy}
\newcommand{\ds}{\displaystyle}
\newcommand{\RR}{\mathbb{R}}
\newcommand{\qed}{\hfill \ensuremath{\square}}
\newcommand{\pd}[2]{\frac {\p #1}{\p #2}}
\newcommand{\Om}{\Omega}
\newcommand{\p}{\partial}
\newcommand{\Z}{\mathbb{Z}}
\newcommand{\N}{\mathbb{N}}

\newcommand{\Ncal}{\mathcal{N}}

\newcommand{\Dcal}{\mathcal{D}}

\newcommand{\ra}{\rangle}
\newcommand{\la}{\langle}
\newcommand{\proof}{\vspace{1ex}\noindent{\em Proof}. \ }
\makeatletter
{\newcommand{\rroman}[1]{\@Roman{#1}}
\makeatother
\newcommand{\ha}{\hat{a}}
\newcommand{\hb}{\hat{b}}

\newtheorem{thm}{Theorem}[section]

\newtheorem{lem}[thm]{Lemma}

\newtheorem{rem}[thm]{Remark}

\def\nm{\noalign{\medskip}}
\baselineskip 18pt
\bibliographystyle{plain}

\numberwithin{equation}{section}
\begin{document}

\title{Reconstructing Small Perturbations of Scatterers from Electric or Acoustic Far-Field Measurements}

\author{ Mikyoung Lim \thanks{\footnotesize Centre de
Math\'{e}matiques Appliqu\'{e}es, Ecole Polytechnique, 91128
Palaiseau Cedex, France (mklim@cmapx.polytechnique.fr,
louati@cmapx.polytechnique.fr). } \and Kaouthar Louati
\footnotemark[1] \and Habib Zribi \thanks{College of Electronics
and Information, Kyung Hee University, Korea (zribi@khu.ac.kr).}}

\maketitle

\begin{abstract} In this paper we consider the problem of
determining the boundary perturbations of an object from far-field
electric or acoustic measurements. Assuming that the unknown
scatterer boundary is a small perturbation of a circle, we develop a
linearized relation between the far-field data and the shape of
the object. This relation is used to find the Fourier coefficients
of the perturbation of the shape.
\end{abstract}

\noindent {\footnotesize {\bf Mathematics subject classification
(MSC2000):} 35R30}


\section{Introduction}

The field of inverse shape problems has been an active research
area for several decades. Several related scalar problems belong
to this field: electric and acoustic scattering form two large
classes. In direct problems one wants to calculate the field
outside a given object. In two common situations, one knows either
the values of the field on the object (the Dirichlet problem), or
the values of the normal derivative of the field on the boundary
(the Neumann problem). Inverse shape problems involve
reconstructing the object shape from measurements of the electric
or acoustic field. Differently from Direct problems which are
usually well posed, inverse problems are ill posed: the solution
has an unstable dependence on the input data.

The formulation of the electric scattering problem is based on the
quasi-static approximation and the related Laplace equation  for
the electric scalar potential. When a perfect conductor is exposed
to extremely low-frequency electric fields, the problem is
equivalent to the Dirichlet boundary value problem for the Laplace
operator.

The sound-soft acoustic scattering problem is characterized by the
condition that the total field vanishes on the boundary of the
scatterer. Thus, acoustic scattering is equivalent to the
Dirichlet boundary value problem for the Helmholtz operator, with
the scattered field equal to the negative of the known incident
field.

These two problems are frequently solved by methods of potential
theory. The single- and double-layer potentials relate a charge
density on the object boundary to the limiting values of
the field and its normal derivative. The resulting integral
equations are then solved in an appropriate function space, a
common choice being the Lebesgue space $L^2$.

In this paper, assuming that the unknown object boundary is  a
small perturbation of a unit circle, we develop for both electric
and acoustic problems a linearized relation between the the
far-field data and the shape of the scatterer. Under this purpose, we
investigate the Dirichlet boundary value problem outside the object entering
the Dirichlet data as parameters and the shape of the object as
variables.

The linearized relation between the far-field data and the object shape is used to find the Fourier
coefficients of the boundary perturbation of the object. Suppose that the angular
oscillations in the perturbation are less than $1/n$. In order to
detect that perturbation, it turns out that one needs to use the
first $n$ eigenvectors of the Dirichlet-to-Neumann operator
corresponding to the unperturbed shape as the Dirichlet boundary
data. We may think that this result is quite general. When the
unknown object is a $C^2$-perturbation of a disk, we obtain
asymptotic formulae for the Dirichlet-to-Neumann operator in terms
of the small perturbations of the object shape, and it
is worth mentioning the expansions of Dirichlet-to-Neumann
operators for rough non-periodic surfaces \cite{milder, RM} and
for periodic interfaces \cite{DF}.

Our approach relies on asymptotic expansions of the far-field data
with respect to the perturbations in the boundary, in much the
same spirit as the recent work \cite{AKLZ} and the text
\cite{bookAK}. We consider only the two-dimensional case, the
extension to three dimensions being obvious. In connection with
our work, we should also mention the paper by Kaup and Santosa
\cite{kaup} on detecting corrosion from steady-state voltage
boundary perturbations and the work by Tolmasky and Wiegmann
\cite{tol} on the reconstruction of small perturbations of an
interface for the inverse conductivity problem.

We deal with electric problems in section 2 and 3, and acoustic
problems in 4 and 5.

\section{Formulation of the Electric Problem}
We consider the reconstructing problem of the perfect conductor $D_\epsilon$ which is the small perturbation
of the unit disk $D$ described by a Lipschitz function $f$ and a small scale factor $\epsilon$, that is
\begin{equation}\label{Depsilon}
 \p{D_\epsilon}(=\p{D}+\epsilon f e_\theta):= \Bigr\{ (1+\epsilon f(\theta))e_{\theta} ~,~~\theta \in [0 ,2\pi] \Bigr\},
 \end{equation}
where $e_{\theta}=(\cos \theta, \sin \theta)$.
\subsection{Electric Scattering problem}
If we apply an initial potential $v^i$ to $\RR^2$ which is homogeneous except for the perfect conductor $D_\epsilon$,
then the derived electric potential $v$ is given by $v=v^i+v^s$, where $v^s$ is the solution to
\begin{equation}\label{usca}
\quad \left\{
\begin{array}{ll}
\ds \Delta v^s = 0,  \quad & \ds \text{in }  \RR^2 \setminus \overline{D}_\epsilon,\\
\nm \ds v^s=-v^i+C(\mbox{constant}), \quad & \ds \text{on}\quad \p{D_\epsilon}, \\
\nm \ds v^s(x)\rightarrow 0,\quad &\mbox{ as }|x|\rightarrow\infty.
\end{array}
\right.
\end{equation}

We denote $v_0^s$ as the perturbation of electric potential due to
the conductor $D$, i.e.,
\begin{equation}\label{usor}
\quad \left\{
\begin{array}{ll}
\ds \Delta v_0^s = 0,  \quad & \ds \text{in }  \RR^2 \setminus \overline{D},\\
\nm \ds v_0^s =-v^i+C(\mbox{constant}), \quad & \ds \text{on}\quad \p{D}, \\
\nm \ds v_0^s(x)\rightarrow 0,\quad &\mbox{ as }|x|\rightarrow\infty,
\end{array}
\right.
\end{equation}
By obtaining a linearized relation between $\epsilon f$ and the
$\epsilon$-order term of $(v^s-v_0^s)(r,\theta)$ as
$r\rightarrow\infty$, we try to recover $D_\epsilon$.

It follows from the Taylor series expansion of $v^i$ near $D$ that
\begin{equation}\label{taylor e}
v^i(1+\epsilon f(\theta),\theta)=v^i(1,\theta)+\epsilon
f\p_rv^i(1,\theta)+O(\epsilon^2).
\end{equation}
Here we used the polar coordinates
$x(r,\theta)=(r\cos\theta,r\sin\theta)$. We investigate $v^s$ by
considering two exterior boundary value problems,
one with Dirichlet value $(-v^i(1,\theta))$ and the other with
$(-\epsilon f\p_rv^i(1,\theta))$ on $\p D_\epsilon$. To do that, we
formulate the fixed boundary value problem.

\subsection{Fixed Dirichlet boundary value problem}
When the boundary value is
prescribed as the $2\pi$-periodic function $\Psi$, the voltage potential outside the conductor $D_\epsilon$ is given
by the harmonic function $u$ which satisfies the following:
\begin{equation}\label{def2}
\quad \left\{
\begin{array}{ll}
\ds \Delta u = 0, \quad & \ds \text{in }  \RR^2 \setminus \overline{D}_\epsilon,\\
\nm \ds u (1+\epsilon f(\theta),\theta) =\Psi(\theta) ,\quad & \mbox{for }\ds\theta\in[0,2\pi],\\
\nm \ds u(x)=\mbox{constant}+O({1}/{|x|}),\quad &\mbox{as
}|x|\rightarrow+\infty.
\end{array}
\right.
\end{equation}

We let $u_0$ be the voltage potential outside the unit disk $D$ with the fixed Dirichlet data $\Psi$ on the boundary, i.e.,
\begin{equation}\label{def1}
\quad \left\{
\begin{array}{ll}
\ds \Delta u_0 = 0,  \quad & \ds \text{in }  \RR^2 \setminus \overline{D},\\
\nm \ds u_0(1,\theta) =\Psi(\theta), \quad & \ds \mbox{for }\theta\in[0,2\pi], \\
\nm \ds u_0(x)=\mbox{constant}+O({1}/{|x|}),\quad &\mbox{as
}|x|\rightarrow+\infty.\end{array} \right.
\end{equation}

We obtain the linearized relation between the boundary interface of the conductor $D_\epsilon$ and $(u-u_0)$ at infinity,
especially when $\Psi$ is given by a $C^4$-function or a Lipschitz function.

\section{Electric far-field formula  and Inversion algorithm}
\subsection{Linearized relation for the Dirichlet problem}
We start by explaining the main idea to obtain the linearized relation.

To derive the asymptotic expansion of the solution $u$ to \eqref{def2} with the given boundary data $\Psi$,
we apply the field expansion method (F.E) (see \cite{DF}). Firstly, we expand  $u$
 in powers of $\epsilon$, i.e.,
\begin{equation}\label{fem}
u(r, \theta)=\sum_{n=0}^{+\infty}u_n(r ,\theta)\epsilon^n.
\end{equation}
Now expanding in terms of $r$ and evaluating \eqref{fem} at
$r=1+\epsilon f$, we obtain that
$$u(1+\epsilon f(\theta),\theta)=u_0(1,\theta)+\epsilon\Bigr(u_1(1,\theta)+\p_ru_0(1,\theta)f(\theta)\Bigr)+O(\epsilon^2).$$

Since $u(1+\epsilon f(\theta),\theta)$ and $u_0(1,\theta)$ have
the same value, $u_1$ can be considered as the decaying harmonic
function which satisfies
\begin{equation}\label{u1field}
u_1(1,\theta)=-\Ncal_0(\Psi)(\theta)f(\theta),
\end{equation}
where $$\Ncal_0(\Psi)(\theta)=\p_ru_0(1,\theta).$$ In other words,
$\Ncal_0$ is the {Dirichlet-to-Neumann operator} of $D$, and it
can be expressed as
\begin{equation}\label{Ncal0}
 \Ncal_0(\Psi)(\theta)=-\sum_{n=1}^{+\infty}\Bigr[n\hat{a}_n(\Psi)\cos n\theta+n\hat{b}_n(\Psi)\sin n\theta \Bigr],
 \end{equation}
 where $\hat{a}_n(\Psi)$ and $\hat{b}_n(\Psi)$ are the fourier coefficients, that is
 $$\hat{a}_n(\Psi)=\frac{1}{\pi}\int_0^{2\pi}\Psi(\theta)\cos (n\theta)  d\theta,\
\hat{b}_n(\Psi)=\frac{1}{\pi}\int_0^{2\pi}\Psi(\theta)\sin (n\theta) d\theta.$$

From \eqref{u1field} and the expansion of the harmonic function
outside a disk, we have
\begin{equation*}
u_1(r,\theta)=-\sum_{n=0}^{+\infty} \frac{1}{r^n}\Bigr[\hat{a}_n(\Ncal_0(\Psi)
f)\cos n\theta+\hat{b}_n(\Ncal_0(\Psi)f)\sin n\theta\Bigr],\quad r\geq1.
\end{equation*}
Therefore,
$$(u-u_0)(r,\theta)\sim -\frac{\epsilon}{r}\Bigr[\hat{a}_1(\Ncal_0(\Psi)
f)\cos\theta+\hat{b}_1(\Ncal_0(\Psi)f)\sin\theta\Bigr]+C,\quad r\gg1,$$
where $C$ is a constant. More precisely, we have the following theorem and give the proof in Subsection \ref{proofs}.
\begin{thm}\label{thm1}
For a $2\pi$-periodic function $\Psi$, we let $u$ and $u_0$ be the solution to \eqref{def2} and \eqref{def1}, respectively.
\begin{enumerate}
\item Let $\Psi\in C^4([0,2\pi])$.  For $r\gg1$, we have
\begin{equation}\label{uasympC}
(u-u_0)(r,\theta)=-\frac{\epsilon}{r}\Bigr[\hat{a}_1(\Ncal_0(\Psi)f
)\cos\theta+\hat{b}_1(\Ncal_0(\Psi)f )\sin\theta\Bigr]+C
+O({\epsilon^\frac{3}{2}}/{r}+{\epsilon}/{r^2}),
\end{equation}
where $C$ is a constant, and $O(\epsilon^\frac{3}{2}/r+\epsilon
/r^2)$ depends on the Lipschitz constant of $f$ and
$\|\Psi\|_{C^4}$.
\item For a Lipschitz function $\Psi$, we have
that
\begin{equation}\label{uasympL}
(u-u_0)(r,\theta)=C+O(\epsilon^\frac{1}{2}/r+1 /r^2),
\quad\mbox{for }r\gg 1,
\end{equation}
where $C$ is a constant, and $O(\epsilon^\frac{1}{2}/r+1/r^2)$
depends on the Lipschitz constant of $f$ and $\Psi$.
\end{enumerate}
\end{thm}

\begin{rem}
For the case of $C^2$-perturbation of the interface, i.e., $f\in
C^2([0,2\pi])$, the error term of \eqref{uasympC} and
\eqref{uasympL} can be replaced by $O(\epsilon^2/r+\epsilon/r^2)$ and $O(\epsilon/r+1/r^2)$.
\end{rem}

In connection with the results for rough non-periodic surfaces
\cite{milder, RM} and for periodic interfaces \cite{DF}, we expand
the Dirichlet-to-Neumann operator $\Ncal_{\epsilon f}$ of
$D_\epsilon$ which is defined by
$$\Ncal_{\epsilon f}(\Psi)(\theta):=\pd{u}{\nu_y}(y),\quad y=(1+\epsilon f(\theta))e_{\theta},$$
where $\nu_y$ is the outward unit normal vector to $D_\epsilon$.

Note that $\nu_y$ is given by
\begin{equation}\label{normal_p}
\nu_y=\frac{N_\theta}{|N_\theta|},
\end{equation}
where
\begin{equation*}
N_\theta=(1+\epsilon f(\theta))e_\theta - \epsilon
\dot{f}\tau_{\theta},\ \tau_{\theta}=(-\sin \theta, \cos \theta).
\end{equation*}
Here $\dot{f}$ is the derivative of $f$ with respect to $\theta$.
From the fact that
\begin{equation}\label{normal_asymp}
\frac{1}{|N_\theta|}=1-\epsilon
f+O(\epsilon^2),
\end{equation} it follows that
\begin{align*}
\Ncal_{\epsilon f}(\Psi)(\theta)&= (1-\epsilon f)\la \nabla u,N_\theta\ra+O(\epsilon^2)\\
&= (1-\epsilon f)\Bigr[(1+\epsilon f)\pd
{u}{r}\Bigr|_{r=1+\epsilon f}-
\frac{{\epsilon \dot{f}}}{1+\epsilon f}\pd {u}{\theta}\Bigr|_{r=1+\epsilon f}\Bigr]+O(\epsilon ^2)\\
&= \pd {u}{r}\Bigr|_{r=1+\epsilon f(\theta)}-\epsilon \dot{f}\pd
{u}{\theta}\Bigr|_{r=1+\epsilon f(\theta)}+O(\epsilon ^2).
\end{align*}
Applying \eqref{fem}, we obtain
\begin{align}\label{asym}
\Ncal_{\epsilon f}(\Psi)(\theta)&\sim\p_r u_0 (1 ,\theta) +\epsilon
\Big(\p_r u_1 (1 ,\theta)+\p_r^2  u_0 (1
,\theta)f(\theta)-\p_{\theta}  u_0 (1
,\theta)\dot{f}(\theta)\Bigr).\nonumber
\end{align}

Defining an operator $\Dcal_0$ by
\begin{equation}\label{Dcal0}
\Dcal_0 (\Psi)(\theta):=-\sum_{n=1}^{+\infty}\Bigr[(n+1)
\hat{a}_n(\Psi)\cos n\theta+(n+1)\hat{b}_n(\Psi)\sin n\theta\Bigr],
\end{equation}
we have
\begin{equation}\label{ur2}
\p_r^2  u_0 (1 ,\theta)=\Dcal_0 \Ncal_0 (\Psi)(\theta).
\end{equation}

\begin{lem}\label{Nepsilon}
For $f\in C^2([0,2\pi])$ and $\Psi \in \mathcal{C}^4([0 , 2\pi])$,
we have
\begin{equation*}
\Ncal_{\epsilon f}(\Psi) =\Ncal_0 (\Psi) + \epsilon \Ncal_f^1
(\Psi)+O(\epsilon^\frac{3}{2}),
\end{equation*}
where
\begin{equation}\label{develop}
\Ncal_f^1 (\Psi)= \Dcal_0 \Ncal_0(\Psi) f -\Ncal_0(\Ncal_0(\Psi)
f) -\dot{f}\dot{\Psi}.
\end{equation}
\end{lem}
We give the proof in Subsection \ref{proofs}.

\subsection{Algorithm for the Inverse Shape Problem}
For an entire harmonic function $v^i$, we let $v^s$ and $v^s_0$ be the solution to \eqref{usca} and \eqref{usor}, respectively.
The Dirichlet values of the solutions are given by
\begin{align*}
v_0^s|_{\p D}&=-v^i(1,\theta)+\mbox{constant},\\
v^s|_{\p D_\epsilon}&=-v^i(1,\theta)-\epsilon
\p_rv^i(1,\theta)f(\theta)+\mbox{constant}+O(\epsilon^2).
\end{align*}

Note that $$\p_rv^i(1,\theta)=-\Ncal_0(v^i|_{\p D}).$$ Here we have the minus sign on the right hand side because $N_0$
is the Dirichlet-to-Neumann operator for the exterior harmonic functions.
Applying \eqref{uasympC} and \eqref{uasympL} with letting $\Psi=-v^i|_{\p D}$ and $\Psi=\Ncal_0(v^i|_{\p D})f$, respectively,
we obtain for $r\gg1$ that
\begin{equation}\label{eqn:farfield}
(v^s-v^s_0)(r,\theta)\sim 2\frac{\epsilon}{r}\Bigr[\hat{a}_1\Bigr(\Ncal_0(v^i|_{\p D})
f\Bigr)\cos\theta+\hat{b}_1\Bigr(\Ncal_0(v^i|_{\p D}) f\Bigr)\sin\theta\Bigr].
\end{equation}

Now define entire harmonic functions $v^{n,i}$ and $w^{n,i}$, for
$n\in \N$, by
 $$v^{n,i}(r,\theta)=-\frac{1}{n}r^n\sin n\theta,\ w^{n,i}(r,\theta)=\frac{1}{n}r^n\cos n\theta.$$
Let $v^{n,s}$ and $w^{n,s}$ be the solution to \eqref{usca} with
the initial potential $v^{n,i}$ and $w^{n,i}$, respectively. In
the same way, define $v_0^{n,s}$ and $w_0^{n,s}$ as the solution to
\eqref{usor}. Let
$$c_1(v^{n,i}):=\frac{1}{\epsilon}\hat{a}_1\Bigr(r\cdot(v^{n,s}-v_0^{n,s})\Bigr),\
d_1(v^{n,i}):=\frac{1}{\epsilon}\hat{b}_1\Bigr(r\cdot(v^{n,s}-v_0^{n,s})\Bigr),$$
$$c_1(w^{n,i}):=\frac{1}{\epsilon}\hat{a}_1\Bigr(r\cdot(w^{n,s}-w_0^{n,s})\Bigr),\
d_1(w^{n,i}):=\frac{1}{\epsilon}\hat{b}_1\Bigr(r\cdot(w^{n,s}-w_0^{n,s})\Bigr).$$

From \eqref{eqn:farfield}, it follows that
\begin{align*}
c_1(v^{n,i})\sim 2 \hat{a}_1\Bigr(f\cdot\Ncal_0(-\frac{1}{n}\sin
n\theta)\Bigr)=2\hat{a}_1\Bigr(f\cdot\sin n\theta\Bigr),\\
d_1(v^{n,i})\sim  2\hat{b}_1\Bigr(f\cdot\Ncal_0(-\frac{1}{n}\sin
n\theta)\Bigr)=2\hat{b}_1\Bigr(f\cdot\sin n\theta\Bigr).
\end{align*}
By the same way, we obtain
\begin{align*}
c_1(w^{n,i})\sim -2\hat{a}_1\Bigr(f\cdot\cos n\theta\Bigr),\\
d_1(w^{n,i})\sim  -2\hat{b}_1\Bigr(f\cdot\cos n\theta\Bigr).
\end{align*}
Thus we obtain that
\begin{align*}
c_1(v^{n,i})\pm d_1(w^{n,i})&=\frac{1}{\pi}\int_0^{2\pi}2f(\theta)(\sin n\theta\cos\theta\mp\cos n\theta\sin\theta)\ d\theta\\
&=\frac{2}{\pi}\int_0^{2\pi}f(\theta)\sin (n\mp1)\theta\ d\theta
=2\hat{b}_{n\mp1}(f),
\end{align*}
\begin{align*}
\pm d_1(v^{n,i})-c_1(w^{n,i})&=\frac{1}{\pi}\int_0^{2\pi}2f(\theta)(\pm\sin n\theta\sin\theta+\cos n\theta\cos\theta)\ d\theta\\
&=\frac{2}{\pi}\int_0^{2\pi}f(\theta)\cos (n\mp1)\theta\ d\theta=2\hat{a}_{n\mp1}(f).
\end{align*}

Therefore, we arrive at
$$\hb_{n-1}(f)=\frac{c_1(v^{n,i})+d_1(w^{n,i})}{2},\qquad \hb_{n+1}(f)=\frac{c_1(v^{n,i})-d_1(w^{n,i})}{2},$$
and
$$\ha_{n-1}(f)=\frac{d_1(v^{n,i})-c_1(w^{n,i})}{2},\qquad \ha_{n+1}(f)=\frac{-d_1(v^{n,i})-c_1(w^{n,i})}{2}, \quad n\geq 1.$$

This simple calculation shows that in order to detect a
perturbation that has oscillations of order $1/n$, one needs to
use the first $n$ eigenvectors ($e^{il \theta}, l=1,\ldots, n,$)
of the Dirichlet-to-Neumann operator $\Ncal_0$ as Dirichlet
boundary data. This is a relatively simple but quite deep
observation. We conjecture that this result holds for general
domains. Another observation is that our asymptotic formula is in
fact a low-frequency expansion which holds for fixed $n$ as
$\epsilon$ goes to zero. It would be interesting to derive an
expansion which is valid for high-frequencies, not just for finite
$n$.

\subsection{Proofs of Theorem \ref{thm1} and Lemma \ref{Nepsilon}}\label{proofs}
We modify $u_0$ and $u_1$ to the solutions $u_0^{\epsilon M}$ and
$u_1^{\epsilon M}$ of
\begin{equation*}
\quad \left\{
\begin{array}{ll}
\ds \Delta u_0^{\epsilon M} = 0, \quad & \ds \text{in }  \RR^2 \setminus \overline{B(1-\epsilon M,0)},\\
\nm \ds u_0^{\epsilon M} (1-\epsilon M,\theta)= \Psi(\theta),\quad & \mbox{for }\ds\theta\in[0,2\pi]\\
\nm \ds u_0(r,\theta)=\mbox{constant}+O({1}/r),\quad &\mbox{as
}r\rightarrow+\infty,
\end{array}
\right.
\end{equation*}
and
\begin{equation*}
\quad \left\{
\begin{array}{ll}
\ds \Delta u_1^{\epsilon M} = 0, \quad & \ds \text{in }  \RR^2 \setminus \overline{B(1-\epsilon M,0)},\\
\nm \ds u_1^{\epsilon M} (1-\epsilon M,\theta)= -[f(\theta)+M]{\p_ru_0(1,\theta)},\quad & \mbox{for }\ds\theta\in[0,2\pi]\\
\nm \ds u_1(r,\theta)=\mbox{constant}+O({1}/r),\quad &\mbox{as
}r\rightarrow+\infty,
\end{array}
\right.
\end{equation*}
where
\begin{equation}\label{fmax}
M:=\max(\|f\|_{L^\infty},\|\dot{f}\|_{L^\infty},1).
\end{equation}

From the fourier expansion of $\Psi$, we obtain
\begin{equation}\label{u0exp}
u_0^{\epsilon M}(r,\theta)
=\sum_{n=0}^{+\infty}\Bigr(\frac{1-\epsilon M}{r}\Bigr)^n
\Bigr[\hat{a}_n(\Psi)\cos n\theta+\hat{b}_n(\Psi)\sin n\theta\Bigr],\quad\mbox{for }r\geq 1-\epsilon M,
\end{equation}
and
\begin{equation}\label{u1epM}
u_1^{\epsilon M}(r,\theta)=-\sum_{n=0}^{+\infty}\Bigr(\frac{1-\epsilon M}{r}\Bigr)^n
\Bigr[\hat{a}_n([f+M]\Ncal_0(\Psi))\cos
n\theta+\hat{b}_n([f+M]\Ncal_0(\Psi))\sin n\theta\Bigr].
\end{equation}

The following is the key lemma to obtain the asymptotic expansion of $(u-u_0)$.
\begin{lem}\label{2asym}
For a $2\pi$-periodic function $\Psi$, we let $u$ and $u_0$ be the solution to \eqref{def2} and \eqref{def1}, respectively.
\begin{enumerate}
\item For $\Psi\in C^4([0,2\pi])$, we have the following
asymptotic expansion holds uniformly on $\p D_\epsilon$:
\begin{equation}\label{lemma:e_linearized}
u=u_0^{\epsilon M}+\epsilon u_1^{\epsilon
M}+C+O(\epsilon^\frac{3}{2}),
\end{equation}
where $C$ is a constant, and $O(\epsilon^\frac{3}{2})$ depends on
the Lipschitz constant of $f$ and $\|\Psi\|_{C^4}$. \item For a
Lipschitz function $\Psi$, we have the following asymptotic
expansion holds uniformly on $\p D_\epsilon$:
\begin{equation}\label{liplinear}
u=u_0^{\epsilon M}+O(\epsilon^\frac{1}{2}),
\end{equation}
where $O(\epsilon^\frac{1}{2})$ depends on the Lipschitz constant of $f$ and $\Psi$.
\end{enumerate}
\end{lem}
\proof Note that the Dirichlet value of $u$ on $\p D_\epsilon $ is $\Psi$. Using \eqref{u0exp} and \eqref{u1epM}, we obtain
\begin{align*}
&(u-u_0^{\epsilon M}-\epsilon u_1^{\epsilon
M})(1+\epsilon f,\theta) \nonumber\\
&=C+\sum_{n=1}^{+\infty}\Bigr[1-\Bigr(\frac{1-\epsilon M}{1+\epsilon
f}\Bigr)^n-\epsilon n(M+f)\Bigr] \Bigr(\hat{a}_n(\Psi)\cos
n\theta+\hat{b}_n(\Psi)\sin
n\theta\Bigr)\\\nonumber
&\quad+\epsilon\sum_{n=1}^{+\infty}\Bigr[\Bigr(\frac{1-\epsilon
M}{1+\epsilon f}\Bigr)^n-1\Bigr]
\Bigr(\hat{a}_n([f+M]\Ncal_0(\Psi))\cos n\theta+\hat{b}_n([f+M]\Ncal_0(\Psi))\sin n\theta\Bigr)\\
&=:C+I+II,\nonumber
\end{align*}
where $C$ is a constant.
\smallskip

Note that
\begin{align}
|1-(1-t)^n-nt|&\leq n^2t^2,\label{poly1}\\\label{poly2}
|1-(1-t)^n|&\leq\max\{1,\ 2nt\}.
\end{align}
For $\Psi\in C^4([0,2\pi])$, we have
\begin{equation}\label{fourmax}
|\hat{a}_n(\Psi)|,|\hat{b}_n(\Psi)|\leq C\frac{\|\Psi\|_{C^4}}{n^4},\qquad\mbox{for each }n\in\N,
\end{equation}
and from \eqref{poly1}, it follows that
$$I=O(\epsilon^2).$$
Now, applying Cauchy-Schwarz inequality, we obtain
\begin{align}
|II|^2&\leq\epsilon^2\sum_{n=1}^{+\infty}\frac{1}{n^2}\nonumber
\Bigr[\Bigr(\frac{1-\epsilon M}{1+\epsilon f}\Bigr)^n-1\Bigr]^2\\
&\qquad\times\sum_{n=1}^{+\infty}n^2\nonumber
\Bigr[\hat{a}_n([f+M]\Ncal_0(\Psi))\cos n\theta+\hat{b}_n([f+M]\Ncal_0(\Psi))\sin n\theta\Bigr]^2\\
&\leq \epsilon^2\Bigr\|\frac{d}{d\theta}\Bigr([f+M]\Ncal_0(\Psi)\Bigr)\Bigr\|^2_{L^2([0,2\pi])} \nonumber
\sum_{n=1}^{+\infty}\frac{1}{n^2}
\Bigr[1-\Bigr(\frac{1-\epsilon M}{1+\epsilon f}\Bigr)^n\Bigr]^2.
\end{align}
From \eqref{poly2}, it follows
\begin{equation}\label{sum2epsilon}
\sum_{n=1}^{+\infty}\frac{1}{n^2}
\Bigr[1-\Bigr(\frac{1-\epsilon M}{1+\epsilon f}\Bigr)^n\Bigr]^2
=\sum_{n\leq1/\epsilon}\frac{1}{n^2}(C\epsilon n)^2
+\sum_{n>1/\epsilon}\frac{1}{n^2}\leq C \epsilon.
\end{equation}
Therefore we have
\begin{equation}\label{gint}
|II|^2\leq C \epsilon^3,
\end{equation}
where $C$ depends on the Lipshitz constant of $f$ and $\|\Psi\|_{C^4}$.
\smallskip

When $\Psi$ is a Lipschtz function, from \eqref{sum2epsilon} we have
\begin{align*}
\Bigr|(u-u_0^{\epsilon M})(1+\epsilon f,\theta)\Bigr|
&=\left|\sum_{n=1}^{+\infty}\Bigr[1-\Bigr(\frac{1-\epsilon M}{1+\epsilon f}\Bigr)^n\Bigr]
\Bigr(\hat{a}_n(\Psi)\cos n\theta+\hat{b}_n(\Psi)\sin n\theta\Bigr)\right|\\
&\leq \|\dot{\Psi}\|_{L^2([0,2\pi])}  \Bigr(\sum_{n=1}^{+\infty}\frac{1}{n^2}
\Bigr[1-\Bigr(\frac{1-\epsilon M}{1+\epsilon f}\Bigr)^n\Bigr]^2\Bigr)^\frac{1}{2}\\\nonumber
&\leq C\epsilon^\frac{1}{2}
\end{align*}
\qed

{\bf{Proof of Theorem \ref{thm1}}}
 For $\Psi\in C^4([0,2\pi])$,
from \eqref{lemma:e_linearized} and the decaying condition of $u$,
$u_0^{\epsilon M}$ and $u_1^{\epsilon M}$ at infinity,
$$u(r,\theta)=(u_0^{\epsilon M}+\epsilon u_1^{\epsilon M})(r,\theta)+\mbox{constant}+O(\epsilon^\frac{3}{2}),\quad r\gg 1.$$
Let $\Om$ be a ball containing $D_\epsilon$, then from the
invertibility of the Double layer potential in $L^2_0(\p\Om)$, it
follows that
$$u(r,\theta)=(u_0^{\epsilon M}+\epsilon u_1^{\epsilon M})(r,\theta)+\mbox{constant}+O(\epsilon^\frac{3}{2}/r),\quad r\gg 1.$$

We calculate that
\begin{align*}
&(u_0^{\epsilon M}+\epsilon u_1^{\epsilon M}-u_0)(r,\theta)\\
&=\sum_{n=0}^{+\infty}\frac{(1-\epsilon
M)^n-1}{r^n}\Bigr[\hat{a}_n(\Psi)\cos n\theta+\hat{b}_n(\Psi)\sin
n\theta\Bigr]\\
&\quad -\epsilon\sum_{n=0}^{+\infty}\Bigr(\frac{1-\epsilon
M}{r}\Bigr)^n \Bigr[\hat{a}_n([f+M]\Ncal_0(\Psi))\cos
n\theta+\hat{b}_n([f+M]\Ncal_0(\Psi))\sin n\theta\Bigr]\\
&=-\frac{\epsilon}{r}\Bigr[\hat{a}_1(\Ncal_0(\Psi)
f)\cos\theta+\hat{b}_1(\Ncal_0(\Psi) f)\sin\theta\Bigr]
+C+O(\frac{\epsilon}{r^2}+\frac{\epsilon^2}{r}),\qquad\mbox{for
}r\gg 1,
\end{align*}
where $C$ is a constant, and $O(\frac{\epsilon}{r^2}+\frac{\epsilon^2}{r})$ depends on the
Lipschtz constant of $f$ and  $\|\Psi\|_{C^4}$. Therefore we prove
\eqref{uasympC}.

 By the same way, we can prove \eqref{uasympL}\qed
\smallskip
\smallskip

{\bf{Proof of Lemma \ref{Nepsilon}}}
Note that
\begin{align}
|\hat{a}_n([f+M]\Ncal_0(\Psi))|,\ |\hat{b}_n([f+M]\Ncal_0(\Psi))|\leq C\frac{1}{n^2},\qquad\mbox{for }n\in\N,\nonumber\\\label{n4fM}
\sum_{n=1}^{+\infty} n^4\Bigr(\hat{a}_n([f+M]\Ncal_0(\Psi))\cos n\theta+\hat{b}_n([f+M]\Ncal_0(\Psi))\sin n\theta\Bigr)^2\leq C,
\end{align}
where $C$ depends on $\|f\|_{C^2}$ and $\|\Psi\|_{C^4}$.

Thus we obtain that
\begin{align}
&\pd{}{\theta}(u-u_0^{\epsilon M}-\epsilon u_1^{\epsilon M})|_{\p D_\epsilon}\nonumber\\\nonumber
&=\sum_{n=1}^{+\infty} n\Bigr[1-\Bigr(\frac{1-\epsilon M}{1+\epsilon f}\Bigr)^n-\epsilon n(M+f)\Bigr]
\Bigr(-\hat{a}_n(\Psi)\sin n\theta+\hat{b}_n(\Psi)\cos n\theta\Bigr)\\\nonumber
&\quad+\epsilon\sum_{n=1}^{+\infty} n\Bigr[\Bigr(\frac{1-\epsilon M}{1+\epsilon f}\Bigr)^n-1\Bigr]
\Bigr(-\hat{a}_n([f+M]\Ncal_0(\Psi))\sin n\theta+\hat{b}_n([f+M]\Ncal_0(\Psi))\cos n\theta\Bigr)\\\nonumber
&\quad+\epsilon\dot{f}(\theta)\sum_{n=1}^{+\infty }n
\Bigr[\Bigr(\frac{1-\epsilon M}{1+\epsilon f}\Bigr)^n\frac{1}{1+\epsilon f}-1\Bigr]
\Bigr(\hat{a}_n(\Psi)\cos n\theta+\hat{b}_n(\Psi)\sin n\theta\Bigr)\\\nonumber
&\quad-\epsilon^2\dot{f}(\theta)\sum_{n=1}^{+\infty} n\Bigr(\frac{1-\epsilon M}{1+\epsilon f}\Bigr)^n\frac{1}{1+\epsilon f}
\Bigr(\hat{a}_n([f+M]\Ncal_0(\Psi))\cos n\theta+\hat{b}_n([f+M]\Ncal_0(\Psi))\sin n\theta\Bigr)\\
&=O(\epsilon^\frac{3}{2}).\label{ptheta}
\end{align}

There exists a constant $C$ which depends on the Lipschitz character of $\p D_\epsilon$, see \cite{bookAK}, such that
\begin{equation*}
\Bigr\|\pd{}{\nu}(u-u_0^{\epsilon M}-\epsilon u_1^{\epsilon M})\Bigr\|_{L^2(\p D_\epsilon)}
\leq C\Bigr\|\pd{}{T}(u-u_0^{\epsilon M}-\epsilon u_1^{\epsilon M})\Bigr\|_{L^2(\p D_\epsilon)},
\end{equation*}
where $T$ is the unit tangent vector on $\p D_\epsilon$. From \eqref{ptheta} and the fact that $\pd{}{\theta}=(1+\epsilon f+O(\epsilon^2))\pd{}{T}$,
it follows that
\begin{equation}\label{NormalTangent}
\Bigr\|\pd{}{\nu}(u-u_0^{\epsilon M}-\epsilon u_1^{\epsilon M})\Bigr\|_{L^2(\p D_\epsilon)}
=O(\epsilon^\frac{3}{2}).
\end{equation}

From \eqref{normal_p}, \eqref{normal_asymp} and \eqref{NormalTangent}, we have
\begin{align}
\Ncal_{\epsilon f}(\Psi)(\theta)
&=\la \nabla (u_0^{\epsilon M}+\epsilon u_1^{\epsilon M}),\  \nu_y\ra+O(\epsilon^\frac{3}{2})\nonumber\\
&= (1-\epsilon f)\la \nabla  (u_0^{\epsilon M}+\epsilon u_1^{\epsilon M}),N_\theta\ra+O(\epsilon^\frac{3}{2})\nonumber\\
&= \pd {}{r} (u_0^{\epsilon M}+\epsilon u_1^{\epsilon M})\Bigr|_{r=1+\epsilon f(\theta)}
-\epsilon \dot{f}\pd {}{\theta} (u_0^{\epsilon M}+\epsilon u_1^{\epsilon M})\Bigr|_{r=1+\epsilon f(\theta)}+O(\epsilon^\frac{3}{2})\nonumber\\
&= \pd {}{r} (u_0^{\epsilon M}+\epsilon u_1^{\epsilon M})\Bigr|_{r=1+\epsilon f(\theta)}
-\epsilon\dot{f}\dot{\Psi}+O(\epsilon^\frac{3}{2}).\label{e:Nepsilon}
\end{align}
We compute
\begin{align*}
&\pd {}{r} (u_0^{\epsilon M}+\epsilon u_1^{\epsilon M})\Bigr|_{r=1+\epsilon f(\theta)}\\
&=\sum_{n=1}^{+\infty}\frac{1}{1+\epsilon f}\Bigr(\frac{1-\epsilon M}{1+\epsilon f}\Bigr)^n
(-n)\Bigr(\hat{a}_n(\Psi)\sin n\theta+\hat{b}_n(\Psi)\cos n\theta\Bigr)\\
&\quad+\epsilon\sum_{n=1}^{+\infty}\frac{1}{1+\epsilon f}\Bigr(\frac{1-\epsilon M}{1+\epsilon f}\Bigr)^n
n\Bigr(\hat{a}_n([f+M]\Ncal_0(\Psi))\cos n\theta+\hat{b}_n([f+M]\Ncal_0(\Psi))\sin n\theta\Bigr).
\end{align*}
Since
$$-n\frac{1}{1+\epsilon f}\Bigr(\frac{1-\epsilon M}{1+\epsilon f}\Bigr)^n=-n+\epsilon n(n+1)f+\epsilon n^2 M+O(\epsilon^2 n^3),$$
$$\epsilon n\frac{1}{1+\epsilon f}\Bigr(\frac{1-\epsilon M}{1+\epsilon f}\Bigr)^n
=\epsilon n+\epsilon n\Bigr[\Bigr(\frac{1-\epsilon M}{1+\epsilon f}\Bigr)^n-1\Bigr]+O(\epsilon^2 n),
$$
we have $$\pd {}{r} (u_0^{\epsilon M}+\epsilon u_1^{\epsilon M})\Bigr|_{r=1+\epsilon f(\theta)}
=\Ncal_0(\Psi)+\epsilon\Bigr(\Dcal_0 \Ncal_0(\Psi) f -\Ncal_0(\Ncal_0(\Psi) f)\Bigr)+O(\epsilon^\frac{3}{2}).$$
From \eqref{e:Nepsilon}, we prove the lemma.

\qed

\section{Formulation of the Acoustic Problem}
Analogously to the Laplacian one, we study the inverse scattering problem of reconstructing a sound-soft obstacle, call it
$D_\epsilon$, whose boundary is the perturbation of the unit circle and is given as \eqref{Depsilon}.
\subsection{Inverse Scattering Problem}
For a incident field $v^i$, we denote $v^s$ and $v_0^s$ as the scattered field from $D_\epsilon$ and $D$, respectively.

In other words, $v^s$ and $v_0^s$ are the solutions to
\begin{equation}\label{def4}
\quad \left\{
\begin{array}{ll}
\ds \Delta v^s + k^2 v^s= 0,\quad & \ds \text{in }  \RR^2 \setminus \overline{ D}_\epsilon ,\\
\nm \ds v^s =-v^i,\quad & \ds \text{on}\quad  \p D_\epsilon, \\
\nm \ds \frac{\p }{\p r}v^s(r ,\theta ) -ik  v^s (r
,\theta)=o({r^{-\frac{1}{2}}}),\quad & \ds r \longrightarrow
+\infty.
\end{array}
\right.
\end{equation}
and
\begin{equation}\label{def3}
\quad \left\{
\begin{array}{ll}
\ds \Delta v_0^s +k^2v_0^s = 0, \quad & \ds \text{in }  \RR^2 \setminus \overline{D} ,\\
\nm \ds v_0^s =-v^i, \quad & \ds \text{on}\quad \p{D},\\
\nm \ds \frac{\p }{\p r}v_0^s(r ,\theta ) -ik  v_0^s (r
,\theta)=o({r^{-\frac{1}{2}}}),\quad & \ds r \longrightarrow
+\infty,
\end{array}
\right.
\end{equation}
Here we used the polar coordinates
$x(r,\theta)=(r\cos\theta,r\sin\theta)$, and the wave number $k$ is given by a positive constant.

\subsection{Fixed Dirichlet boundary value problem}
For a $2\pi$-periodic continuous function $\Psi$, we let $u$ be the solution to the Helmholtz problem
with the prescribed boundary data $\Psi$ on $\p D_\epsilon$,
i.e.,
\begin{equation}\label{udef4}
\quad \left\{
\begin{array}{ll}
\ds \Delta u + k^2 u= 0,\quad & \ds \text{in }  \RR^2 \setminus \overline{D}_\epsilon ,\\
\nm \ds u (1+\epsilon f(\theta),\theta) =\Psi (\theta),\quad & \mbox{for }\ds\theta\in[0,2\pi], \\
\nm \ds \frac{\p }{\p r}u(r ,\theta ) -ik  u (r
,\theta)=o({r^{-\frac{1}{2}}}),\quad & \ds r \longrightarrow
+\infty.
\end{array}
\right.
\end{equation}

The solution $u$ corresponding to the unit disk $D$ satisfies that
\begin{equation}\label{u0def3}
\quad \left\{
\begin{array}{ll}
\ds \Delta u_0 +k^2u_0 = 0, \quad & \ds \text{in }  \RR^2 \setminus \overline{D} ,\\
\nm \ds u_0(1,\theta) =\Psi (\theta),\quad & \mbox{for }\ds\theta\in[0,2\pi],\\
\nm \ds \frac{\p }{\p r}u_0(r ,\theta ) -ik  u_0 (r
,\theta)=o({r^{-\frac{1}{2}}}),\quad & \ds r \longrightarrow
+\infty.
\end{array}
\right.
\end{equation}

We investigate the Far-Field difference between $u$ and $u_0$, especially when $\Psi$ is a $C^4$-function or a Lipschitz function.
\section{Acoustic far-field formula  and inversion algorithm}
\subsection{Asymptotic Far-field expansion for the Dirichlet problem}
We parametrize the unit circle $\p D$ by $\theta\in[0,2\pi]$ and expand $\Psi$ as
$$\Psi(\theta)=\sum_{n\in \mathbb{Z}}\hat{c}_n(\Psi) e^{i n \theta},$$
where $\hat{c}_n(\Psi)$ is the fourier coefficient with respect to $e^{in\theta}$.
By the uniqueness of the exterior Dirichlet problem, it follows that
\begin{equation}\label{u0scattering}
u_0(r,\theta) = \sum_{n \in \Z}\frac{H^{(1)}_{|n|}(k r )}{H^{(1)}_{|n|}(k)}
\ \hat{c}_n(\Psi) e^{i  n \theta}.
\end{equation}

Define the {Dirichlet-to-Neumann operator} $\Ncal_{0}$ with respect to $D$ by
$$\Ncal_{0} : u_0|_{S} \to  \p_r u_0|_{S},$$ then, in a
pseudodifferential fashion, $\Ncal_0$ can be written as follows (see \cite{KS}):
\begin{equation}\label{Ncal0H}
\Ncal_{0}(\Psi)(\theta)=\sum_{n\in \mathbb{Z}}\sigma_{1}(n,k)\hat{c}_n(\Psi) e^{i n \theta},
\end{equation}
where the so-called {discrete symbol} $\sigma_{1}$ is given by
\begin{equation*}
\sigma_{1}(n,k)=k \frac{H^{(1)'}_{|n |}(k)}{H^{(1)}_{|n |} (k)}=-k \frac{H^{(1)}_{|n+1|}(k)}{H^{(1)}_{|n |} (k)}+|n|.
\end{equation*}
Thus, for fixed $k$, we have
\begin{equation}
\sigma_1(n,k)\sim |n|,\quad\mbox{as }|n|\rightarrow\infty.
\end{equation}

By the same way as the electric problem, we
can consider $u_1$ as the solution to \eqref{u0def3} with the boundary value $(-\Ncal_0(\Psi)f)$ on $\p D$ instead of $\Psi$, and
it follows
\begin{align*}
u_1(r,\theta)&=-\sum_{n \in\Z}\frac{H^{(1)}_{|n|}(kr)}{H^{(1)}_{|n|}(k)}\hat{c}_n(\Ncal_0(\Psi)f)e^{in\theta}.
\end{align*}
It is known that, for a fixed $n$, the Hankel function of the first kind satisfies
\begin{equation}\label{hankel} H_{|n|}^{(1)}(x)= \sqrt{\frac {2}{\pi x}}e^{i(x- \frac{\pi}{4}-|n |\frac{\pi}{2})}
+O(|x|^{-1}),\quad x\gg |n|.\end{equation}
We refer to \cite{AS} for more properties of the Hankel function.

Choose $N\in\N$ satisfying that
\begin{equation}\label{bigN}
\sum_{|n|>N}\Bigr|\hat{c}_n(f\Ncal_0(\Psi))\Bigr|=O(\epsilon^\frac{1}{2}),
\end{equation}
then we have the following lemma. More precise proof is given in the Subsection \ref{acoustic:proof}.
\begin{thm}\label{thm:sca}
\begin{enumerate}
\item Let $\Psi\in C^4([0,2\pi])$ and $u$ be the solution to \eqref{udef4}. For $r\gg1$, we have
\begin{equation}\label{acous_u1}
(u-u_0)(r,\theta)= -\epsilon\sqrt\frac{2}{\pi
r}e^{ikr}\sum_{|n|\leq N}
\frac{\hat{c}_n(\Ncal_0(\Psi)f)}{H^{(1)}_{|n|}(k)}e^{-i(\frac{\pi}{4}+\frac{|n|\pi}{2})}e^{in\theta}
+O({\epsilon^\frac{3}{2}}/{\sqrt r}+{\epsilon}/{r}).
\end{equation}
where $N$ is defined by \eqref{bigN}, and $O(\epsilon^\frac{3}{2}/\sqrt r+\epsilon/r)$ depends on the Lipschitz constant of $f$ and $\|\Psi\|_{C^4}$.
\item Let $\Psi$ be a Lipschitz function and $u$ be the solution to \eqref{udef4}. We have that
\begin{equation}\label{acous_u2}
(u-u_0)(r,\theta)=O({\epsilon^\frac{1}{2}}/{\sqrt r}),
\quad\mbox{for }r\gg 1,
\end{equation}
where $O(\epsilon^\frac{1}{2}/\sqrt r)$ depends on the Lipschitz constant of $f$ and $\Psi$.
\end{enumerate}
\end{thm}
\smallskip
\smallskip

We define
\begin{equation}\label{Dcal0H}
\Dcal_0 (\Psi)(\theta):=\sum_{n\in \Z}\sigma_2(n,k)\hat{c}_n(\Psi)
e^{i n \theta},
\end{equation}
with $$\sigma_2(n,k)=k \frac{H^{(1)''}_{|n |}(k)}{H^{(1)'}_{|n |}
(k)}.$$

\begin{lem}\label{Nh_ep}
For $f\in C^2([0,2\pi])$ and $\Psi \in C^4([0
, 2\pi])$, we have that
$$\Ncal_{\epsilon f} (\Psi)=\Ncal_0 (\Psi) + \epsilon \Ncal_f^1 (\Psi)+O(\epsilon^\frac{3}{2}),$$
where
$$\Ncal_f^1 (\Psi)= \Dcal_0 \Ncal_0(\Psi) f -\Ncal_0(\Ncal_0(\Psi) f)-\dot{f}\dot{\Psi}.$$
\end{lem}
\subsection{Algorithm for the Inverse Shape Problem}
Let $v^i$ be the incoming wave, and define $v^s$ and $v_0^s$ as the solution to \eqref{def4} and \eqref{def3}, respectively.
Note that
\begin{equation*}\label{taylor e}
v^i(1+\epsilon f(\theta), \theta)=v^i(1,\theta)+\epsilon
f(\theta)\p_rv^i(1,\theta)+O(\epsilon^2).
\end{equation*}

Applying Theorem \ref{thm:sca} by letting $\Psi=-v^i(1,\theta)$ and $\Psi=-\epsilon
f(\theta)\p_rv^i(1,\theta)$, we have for $r\gg1$ that
\begin{equation*}
(v^s-v_0^s)(r,\theta)\sim \epsilon\sqrt\frac{2}{\pi r}e^{ikr}
\sum_{|n|\leq N}\frac{\hat{c}_n\Bigr(f\Ncal_0(v^i|_{\p D})-f\p_rv^i|_{\p D}\Bigr)}{H^{(1)}_{|n|}(k)}
e^{-i(\frac{\pi}{4}+\frac{n\pi}{2})}e^{in\theta},
\end{equation*}
where $N$ is defined by
$
\sum_{|n|>N}\Bigr|\hat{c}_n\Bigr(f\Ncal_0(v^i|_{\p D})-f\p_rv^i|_{\p D}\Bigr)\Bigr|=O(\epsilon^\frac{1}{2}).
$
This yields to stable reconstruction of
the Fourier coefficients $\hat
c_{n}\Bigr(f\Ncal_0(v^i|_{\p D})-f\p_rv^i|_{\p D}\Bigr)$ for $n$ such
that $H^{(1)}_{|n|}(k)$ is not too big.

\smallskip

Suppose now that $v^i$ satisfies
$$\Ncal_0(v^i|_{\p D})-\p_rv^i|_{\p D}=e^{-i(m-1)\theta},$$
then by measuring $\hat c_{1}(v^s-v_0^s)$, then we can
reconstruct $\hat c_{m}(f)$.

\subsection{Proofs of Theorem \ref{thm:sca} and Lemma \ref{Nh_ep}}\label{acoustic:proof}
We modify $u_0$ and $u_1$ as $u_0^{\epsilon M}$ and
$u_1^{\epsilon M}$ which satisfy
\begin{equation*}
\quad \left\{
\begin{array}{ll}
\ds (\Delta+k^2) u_0^{\epsilon M} = 0, \quad & \ds \text{in }  \RR^2 \setminus \overline{B(1-\epsilon M,0)},\\
\nm \ds u_0^{\epsilon M} (1-\epsilon M,\theta)= \Psi(\theta),\quad & \mbox{for }\ds\theta\in[0,2\pi]\\
\nm \ds \frac{\p }{\p r}u_0^{\epsilon M}(r ,\theta ) -ik  u_0^{\epsilon M} (r
,\theta)=o({r^{-\frac{1}{2}}}),\quad & \ds r \longrightarrow
+\infty.\end{array}
\right.
\end{equation*}
and
\begin{equation*}
\quad \left\{
\begin{array}{ll}
\ds (\Delta+k^2) u_1^{\epsilon M} = 0, \quad & \ds \text{in }  \RR^2 \setminus \overline{B(1-\epsilon M,0)},\\
\nm \ds u_1^{\epsilon M} (1-\epsilon M,\theta)= -[f(\theta)+M]{\Ncal_0(\Psi)},\quad & \mbox{for }\ds\theta\in[0,2\pi],\\
\nm \ds \frac{\p }{\p r}u_1^{\epsilon M}(r ,\theta ) -ik  u_1^{\epsilon M} (r
,\theta)=o({r^{-\frac{1}{2}}}),\quad & \ds r \longrightarrow
+\infty,\end{array}
\right.
\end{equation*}
where $M$ is the constant defined by \eqref{fmax}.
From the fourier expansion of $\Psi$ and the uniqueness of the exterior Dirichlet problem,
$u_0^{\epsilon M}$ and $u_1^{\epsilon M}$ have the expansion as follows:
\begin{equation}\label{Hu0e}
u_0^{\epsilon M}(r,\theta) = \sum_{n \in \Z}\frac{H^{(1)}_{|n|}(k r )}{H^{(1)}_{|n|}(k-\epsilon kM)}
\ \hat{c}_n(\Psi) e^{i  n \theta},
\end{equation}
\begin{equation}\label{Hu1}
u_1^{\epsilon M}(r,\theta)=- \sum_{n \in \Z}\frac{H^{(1)}_{|n|}(k r )}{H^{(1)}_{|n|}(k-\epsilon kM)}
\ \hat{c}_n([f+M]\Ncal_0(\Psi)) e^{i  n \theta}.
\end{equation}

We have the following key lemma to prove Theorem \ref{thm:sca} and Lemma \ref{Nh_ep}.
\begin{lem}\label{acous:2asym}
\begin{enumerate}
\item For $\Psi\in C^4([0,2\pi])$, we have the following
asymptotic expansion for the solution $u$ to \eqref{udef4} holds uniformly on $\p D_\epsilon$:
\begin{equation}\label{linearized}
u=u_0^{\epsilon M}+\epsilon u_1^{\epsilon
M}+O(\epsilon^\frac{3}{2}),
\end{equation}
where $O(\epsilon^\frac{3}{2})$ depends on
the Lipschitz constant of $f$ and $\|\Psi\|_{C^4}$.
 \item For a
Lipschitz function $\Psi$,we have the following asymptotic
expansion for the solution $u$ to \eqref{udef4} holds uniformly on $\p D_\epsilon$:
\begin{equation}\label{liplinear}
u=u_0^{\epsilon M}+O(\epsilon^\frac{1}{2}),
\end{equation}
where $O(\epsilon^\frac{1}{2})$ depends on the Lipschitz constant of $f$ and $\Psi$.
\end{enumerate}
\end{lem}
\proof
From \eqref{Hu0e}, \eqref{Hu1} and the boundary condition of $u$ on $\p D_\epsilon$, we have that
\begin{align*}
(u-u_0^{\epsilon M}-\epsilon u_1^{\epsilon M})(1+\epsilon f(\theta),\theta)
&=\sum_{n \in \Z}\Bigr[1-\frac{H^{(1)}_{|n|}(k +\epsilon kf )}{H^{(1)}_{|n|}(k-\epsilon kM)}+\epsilon [f+M]\sigma_1(n,k)\Bigr]
\ \hat{c}_n(\Psi) e^{i  n \theta}\\
&\quad+\epsilon \sum_{n \in \Z}\Bigr[\frac{H^{(1)}_{|n|}(k+\epsilon kf )}{H^{(1)}_{|n|}(k-\epsilon kM)}-1\Bigr]
\ \hat{c}_n([f+M]\Ncal_0(\Psi)) e^{i  n \theta},
\end{align*}
and
\begin{align*}
(u-u_0^{\epsilon M})(1+\epsilon f(\theta),\theta)
&=\sum_{n \in \Z}\Bigr[1-\frac{H^{(1)}_{|n|}(k +\epsilon kf )}{H^{(1)}_{|n|}(k-\epsilon kM)}\Bigr]
\ \hat{c}_n(\Psi) e^{i  n \theta}.
\end{align*}

Note that
$$\left|H^{(1)}_{|n|}(k+\epsilon t)-H^{(1)}_{|n|}(k)\right|
\leq \epsilon|t|\ \|H^{(1)'}_{|n|}\|_{L^\infty([k-\epsilon |t|,\ k+\epsilon |t|])},$$
$$\left|H^{(1)}_{|n|}(k+\epsilon t)-H^{(1)}_{|n|}(k)-\epsilon t H^{(1)'}_{|n|}(k)\right|
\leq\frac{\epsilon^2t^2}{2}\|H^{(1)''}_{|n|}\|_{L^\infty([k-\epsilon|t|,\ k+\epsilon |t|])}.$$
From the fact that
$$H^{(1)'}_{|n|}(z)=-H^{(1)}_{|n+1|}(z)+\frac{|n|}{z}H^{(1)}_{|n|}(z),$$
we can show that
\begin{align}
1-\frac{H_{|n|}^{(1)}(k+\epsilon kf)}{H_{|n|}^{(1)}(k-\epsilon kM)}
&=O(\epsilon n),\label{Hnbig2}\\
1-\frac{H_{|n|}^{(1)}(k+\epsilon kf)}{H_{|n|}^{(1)}(k-\epsilon kM)}+
\epsilon [f+M]\sigma_1(n,k)&=O(\epsilon^2n^2),\label{Hnbig1}
\end{align}
where $O(\epsilon n)$ and $O(\epsilon^2n^2)$ depend on $M$ and $k$.
Moreover, $|H_{|n|}^{(1)}|(z)$ is decreasing function for $z>0$, and
\begin{equation}\label{Hnbig0}
1-\frac{H_{|n|}^{(1)}(k+\epsilon kf)}{H_{|n|}^{(1)}(k-\epsilon kM)}=O(1).
\end{equation}

Using \eqref{Hnbig2}, \eqref{Hnbig1} and \eqref{Hnbig0}, we can prove the lemma by the same way to prove Lemma \ref{2asym}.
\qed

%
{\bf{Proof of Theorem \ref{thm:sca}}} At first, we assume $\Psi\in
C^4$.
From the solution expression using boundary integral methods (for example, see \cite{CK}), we can show that
\begin{equation}\label{atinfinity}
(u-u_0^{\epsilon M}-\epsilon u_1^{\epsilon M})(r,\theta)=\frac{1}{\sqrt{r}}
\|u-u_0^{\epsilon M}-\epsilon u_1^{\epsilon M}\|_{L^2(\p D_\epsilon)},\quad\mbox{as }r\rightarrow +\infty.
\end{equation}
Therefore $$(u-u_0)(r,\theta)=(u_0^{\epsilon M}+\epsilon u_1^{\epsilon M}-u_0)(r,\theta)+O({\epsilon^\frac{3}{2}}/{\sqrt r}),\quad\mbox{as }r\rightarrow +\infty.$$
From \eqref{u0scattering}, \eqref{Hu0e} and \eqref{Hu1}, we compute that
\begin{align*}
&(u_0^{\epsilon M}+\epsilon u_1^{\epsilon M}-u_0)(r,\theta)\\
&=\sum_{n \in \Z}\frac{H^{(1)}_{|n|}(kr)}{H^{(1)}_{|n|}(k)}\Bigr[\frac{H^{(1)}_{|n|}(k)}{H^{(1)}_{|n|}(k-\epsilon kM)}-1\Bigr]
\ \hat{c}_n(\Psi) e^{i  n \theta}\\
&\quad+\epsilon \sum_{n \in \Z}\frac{H^{(1)}_{|n|}(kr)}{H^{(1)}_{|n|}(k)}
\Bigr[-\frac{H^{(1)}_{|n|}(k)}{H^{(1)}_{|n|}(k-\epsilon kM)}\Bigr]
\ \Bigr(\hat{c}_n(f\Ncal_0(\Psi)) +\hat{c}_n(M\Ncal_0(\Psi))\Bigr)e^{i  n \theta}\\
&=\sum_{n \in \Z}\frac{H^{(1)}_{|n|}(kr)}{H^{(1)}_{|n|}(k)}\Bigr[\frac{H^{(1)}_{|n|}(k)}{H^{(1)}_{|n|}(k-\epsilon kM)}-1-\epsilon M\sigma_1(n,k)\Bigr]
\ \hat{c}_n(\Psi) e^{i  n \theta}\\
&\quad+\epsilon \sum_{n \in \Z}\frac{H^{(1)}_{|n|}(kr)}{H^{(1)}_{|n|}(k)}
\Bigr[-\frac{H^{(1)}_{|n|}(k)}{H^{(1)}_{|n|}(k-\epsilon kM)}\Bigr]
\ \hat{c}_n(f\Ncal_0(\Psi)) e^{i  n \theta}\\
&\quad+\epsilon \sum_{n \in \Z}\frac{H^{(1)}_{|n|}(kr)}{H^{(1)}_{|n|}(k)}
\Bigr[1-\frac{H^{(1)}_{|n|}(k)}{H^{(1)}_{|n|}(k-\epsilon kM)}\Bigr]
\ M\hat{c}_n(\Ncal_0(\Psi)) e^{i  n \theta}\\
&=:I+II+III.
\end{align*}

 From \eqref{Hnbig2} and \eqref{Hnbig1} with replacing $f$ by 0 and \eqref{atinfinity}, it follows
\begin{align}
&I+III=\frac{1}{\sqrt{r}}O(\epsilon^2).\label{Hankel13}
\end{align}
Using \eqref{hankel} and \eqref{atinfinity}, we obtain
\begin{align}
&II=-\epsilon\sqrt\frac{2}{\pi r}e^{ikr}\sum_{|n|\leq N}
\frac{\hat{c}_n(f\Ncal_0(\Psi))}{H^{(1)}_{|n|}(k)}e^{-i(\frac{\pi}{4}+\frac{|n|\pi}{2})}e^{in\theta}+\frac{1}{r}O(\epsilon)
+O(\epsilon^\frac{3}{2}/\sqrt r)
.\label{Hankel2}
\end{align}

When $\Psi$ is a Lipschitz function, we have
$$(u-u_0)(r,\theta)=(u_0^{\epsilon M}-u_0)(r,\theta)+O({\epsilon^\frac{1}{2}}/{\sqrt r})
=O({\epsilon^\frac{1}{2}}/{\sqrt r}),\quad\mbox{as }r\rightarrow +\infty.$$
\qed

\smallskip
{\bf{Proof of Lemma \ref{Nh_ep}}} Note that $\p D$ is a
$C^2$-domain, and using boundary integral methods, we have that (see \cite{CK})
\begin{equation}
\Bigr\|\pd{}{\nu}(u-u_0^{\epsilon M}-\epsilon u_1^{\epsilon
M})\Bigr\|_{C^{0,\alpha}(\p D_\epsilon)} \leq
C\Bigr\|u-u_0^{\epsilon M}-\epsilon u_1^{\epsilon
M}\Bigr\|_{C^{1,\alpha}(\p D_\epsilon)}.
\end{equation}
By the same way as the conductivity case, we can prove the lemma by
calculating $\pd{}{\nu}(u_0^{\epsilon M}+\epsilon u_1^{\epsilon
M})$. \qed

\section*{Acknowledgements} The authors would like to express the gratitude to
\textrm{Professor Habib Ammari} for his kind help on various
points.


\begin{thebibliography}{99}
\bibitem{bookAK} H. Ammari and H. Kang,
{\it Reconstruction of Small Inhomogeneities  from Boundary
Measurements}, Lecture Notes in Mathematics, Volume 1846,
Springer-Verlag, Berlin, 2004.
\bibitem{AKLZ} H. Ammari, H. Kang, M. Lim, and H. Zribi, Conductivity interface problems. Part I: Small perturbations of an
interface, preprint.

\bibitem{AS}M. Abrohmwitz and I. A. Stegun,
\newblock {\it Handbook of Mathematical Functions},
\newblock New York: Dover, 1974.


\bibitem{RM} R. Coifman, M. Goldberg, T. Hrycak, M. Israeli, and V. Rokhlin,
 \newblock An improved operator expansion algorithm for direct and inverse scattering computations,
 \newblock  Waves Random Media 9 (1999), 441-457.

\bibitem{CK} D. Colton and R. Kress,
\newblock {\it Inverse Acoustic and Electromagnetic Scattering Theory},
\newblock 2nd ed. Berlin, Germany: Springer-Verlag, 1998.

\bibitem{kaup} P.G. Kaup and F. Santosa, Nondestructive evaluation
of corrosion damage using electrostatic measurements, J. Nondestr.
Eval. 14 (1995), 127-136.

\bibitem{KS}  M. F. Kondratieva and S. Yu. Sadov,
\newblock Symbol of the Dirichlet-to-Neumann operator in 2D diffraction problems with large wavenumber,
\newblock Day on Diffraction, 2003 Proceedings. International Seminar, 88 - 98.

\bibitem{milder} D.M. Milder, An improved formalism for wave scattering
from rough surfaces, J. Acoust. Soc. Am., 89 (1991), 529–541.
\bibitem{DF} D.P. Nicholls and F. Reitich,
\newblock Analytic continuation of Dirichlet-Neumann operators,
\newblock  Numer. Math. 94 (2003), 107-146.

\bibitem{tol} C.F. Tolmasky and A. Wiegmann, Recovery of small
perturbations of an interface for an elliptic inverse problem via
linearization, Inverse Problems 15 (1999), 465-487.




\end{thebibliography}
\end{document}